\input{amssymb.sty}

\documentclass[11pt]{amsart}
\usepackage{amscd}

\theoremstyle{plain}
\newtheorem{thm}{Theorem}[section]
\newtheorem*{thm*}{Theorem}
\newtheorem{cor}[thm]{Corollary}

\newtheorem*{cor*}{Corollary}
\newtheorem*{conj*}{Conjecture}
\newtheorem*{lemma*}{Lemma}
\newtheorem{lemma}[thm]{Lemma}
\newtheorem*{prop*}{Proposition}
\newtheorem{prop}[thm]{Proposition}

\theoremstyle{definition}

\newtheorem*{defn*}{Definition}
\newtheorem*{rems*}{Remarks}
\newtheorem*{proof*}{Proof}
\newtheorem{prel*}{Preliminaries}
\newtheorem{examples*}{Examples}
\newtheorem{rem}[thm]{Remark}
\newtheorem*{rem*}{Remark}

\newcommand{\C}{\mathbb{C}}

\newcommand{\F}{\widetilde{\mathcal F}}
\newcommand{\M}{\widetilde{M}}

\newcommand{\D}{{\mathbb D}}

\def\ga{\gamma}

\def\eps{\varepsilon}
\def\ka{\kappa}

\def\spec{\hbox{spec}}
\def\Dom{\operatorname{Dom}}

\def\Ker{\operatorname{Ker}\,}

\def\rank{\operatorname{rank}}

\def\sup{\operatorname{sup}}

\def\C{\mathbb C}
\def\R{\mathbb R}

\def\Tr{\operatorname{Tr}}

\def\A{{\mathcal A}} \def\B{{\mathcal B}}
 
\def\F{{\mathcal F}}

\def\tr{\operatorname{tr}}

\def\H{\mathcal H}

\def\D{\mathcal D}

\def\U{{\mathcal U}}
\def\<{\langle}
\def\>{\rangle}

\newcommand{\supp}{\operatorname{supp}}

\newcommand{\nc}{\newcommand}
\nc{\nt}{\newtheorem}
\nc{\gf}[2]{\genfrac{}{}{0pt}{}{#1}{#2}}
\nc{\mb}[1]{{\mbox{$ #1 $}}}
\nc{\real}{{\mathbb R}}
\nc{\comp}{{\mathbb C}}
\nc{\ints}{{\mathbb Z}}
\nc{\Ltoo}{\mb{L^2({\mathbf H})}}
\nc{\rtoo}{\mb{{\mathbf R}^2}}
\nc{\slr}{{\mathbf {SL}}(2,\real)}
\nc{\slz}{{\mathbf {SL}}(2,\ints)}
\nc{\su}{{\mathbf {SU}}(1,1)}
\nc{\so}{{\mathbf {SO}}}
\nc{\hyp}{{\mathbb H}}
\nc{\disc}{{\mathbf D}}
\nc{\torus}{{\mathbb T}}

\nc{\ca}{{\mathcal A}}
\nc{\cag}{{{\mathcal A}^\Gamma}}
\nc{\cg}{{\mathcal G}}
\nc{\chh}{{\mathcal H}}
\nc{\ck}{{\mathcal B}}
\nc{\cm}{{\mathcal M}}
\nc{\cs}{{\mathcal S}}
\nc{\cz}{{\mathcal Z}}
\nc{\sind}{\sigma{\rm -ind}}

\begin{document}

\title[SEMICLASSICAL ASYMPTOTICS]  {SEMICLASSICAL ASYMPTOTICS\\
    AND GAPS IN THE SPECTRA\\
    OF MAGNETIC SCHR\"ODINGER OPERATORS}
\author{V. Mathai}
\address{Department of Mathematics, University of Adelaide, Adelaide
5005, Australia;}
\address{Department of Mathematics, MIT, Cambridge, MA 02139, USA}
\email{vmathai@maths.adelaide.edu.au, \; vmathai@math.mit.edu}
\author{M. Shubin}
\address{Department of Mathematics, Northeastern University, Boston,
MA 02115, USA}
\email{shubin@neu.edu}

\thanks{V.M. acknowledges support from the Clay Mathematical Institute.
M.S. acknowledges partial support from NSF grant DMS-9706038.}


\begin{abstract} In this paper, we study an $L^2$ version of the semiclassical
approximation  of magnetic Schr\"odinger operators with invariant
Morse type potentials on covering spaces of compact manifolds.
In particular, we are  able to establish the existence of an arbitrary large
number of  gaps in the spectrum of these
operators, in the semiclassical limit as the coupling constant $\mu$
goes to zero.
\end{abstract}
\maketitle

\section*{Introduction}
A charged particle constrained to a manifold, moves along geodesic
orbits. However,
in the presence of a magnetic field, the orbit of the charged
particle is no longer
a geodesic, but rather a {\em magnetic geodesic}, that is, a trajectory
determined by solving the Hamiltonian equations given by the symbol
of the magnetic Hamiltonian. In the case of Euclidean space
and constant magnetic field, the magnetic geodesics are circles,
and on the hyperbolic space they are the lines of tangency of the hyperboloids
embedded in the Minkowski space. The corresponding magnetic
Schr\"odinger operators  are sometimes also
called {\em magnetic Hamiltonians} or Landau Hamiltonians. The magnetic
Schr\"odinger operators turn out to possess rich spectral properties
cf. \cite{Bel}, \cite{BrSu},   \cite{Comtet}, \cite{Comtet+H},
\cite{Ma}, that are important in the analysis of the
quantum Hall effect cf. \cite{Bel+E+S},  \cite{CHMM}, \cite{MM},
\cite{Xia}.

We begin by reviewing the construction of the magnetic Hamiltonian on
the universal covering
space $\widetilde M$ of a compact connected manifold
$M$. The fundamental group $\Gamma=\pi_1(M)$ acts on $\widetilde M$ by the deck
transformations, so that $\widetilde M/\Gamma=M$.
  Let $\omega$  be a closed 2-form on $M$ and $B$ be
its lift to $\M$, so
that $B$ is a $\Gamma$-invariant closed 2-form on $\M$. We assume
that $B$ is exact.  Pick
a 1-form $A$ on $\M$ such that $dA = B$. As in geometric
quantization we may regard $A$ as
defining a Hermitian connection $\nabla = d+iA$ on the trivial   line
bundle $\mathcal L$
over $\M$, whose curvature is $iB$. Physically we can think of $A$ as
the electromagnetic
vector potential for a magnetic field $B$. Using the Riemannian
metric the magnetic Laplacian
is given by
$$
H_A = \nabla^*\nabla = (d+iA)^*(d+iA).
$$
Then the Schr\"odinger equation describing the
quantum mechanics of  a single electron of mass $m$ which is confined
to move on the covering space
$\widetilde M$ in the presence of a periodic magnetic field is given by
$$
i \hbar\frac{\partial}{\partial t} \psi = \frac{1}{2m} \left(\hbar 
d+i e A \right)^*
\left(\hbar d+i e A \right)\psi + \mu^{-2} V \psi
$$
where $H_A$ is the magnetic Laplacian, $V$ is a
$\Gamma$-invariant electric potential function, $\hbar$ is Planck's constant,
$e$ is the electric charge of the charge carrier
and $\mu$ is the coupling constant.  In the time
independent framework,  the relevant operator is the magnetic
Schr\"odinger operator
$$
H_{A,V}(\mu) = \mu H_A  + \mu^{-1} V,
$$
where the physical constants are set equal to $1$.
It is the qualitative aspects of
the spectrum of $H_{A,V}(\mu)$ that are relevant to
the study of the  quantum Hall effect. An important feature
of the magnetic
Schr\"odinger operator is that it commutes with magnetic
translations, that is with a projective action of the
fundamental group $\Gamma$. Let $\sigma$ denote the
multiplier or $U(1)$-valued 2-cocycle on $\Gamma$
defining this projective action.
Under the assumption that the Kadison constant of $C^*_r(\Gamma, \sigma)$
is positive, it was proved by Br\"uning and Sunada \cite{BrSu} that
there are only a
finite  number of gaps in the spectrum of the magnetic Schr\"odinger operator
$H_{A,V}$ that lie in any left half-line $(-\infty, \lambda]$. In fact,
under the same hypotheses, they obtain Weyl-type asymptotics of the
number of gaps in the spectrum, as $\lambda \to \infty$. On the other
hand, it was proved in \cite{Ma} that the Kadison constant is positive
under the assumptions that the Baum-Connes
conjecture with coefficients holds for $\Gamma$, and
that $\Gamma$ has finite cohomological dimension and finally that
$\sigma$ defines a rational cohomology class.
The assumptions on $\Gamma$ are satisfied for instance when
$\Gamma$ is a discrete subgroup of $SO(1,n)$ or $SU(1,n)$
or of an amenable Lie group. Our paper on the other hand
proves that $H_{A,V}$ can have arbitrarily large number of gaps,
if $V$ is chosen to be a suitable Morse-type potential. An outstanding open
problem is to construct magnetic Schr\"odinger operators
$H_{A,V}$ that have infinitely many gaps in some half-line
$(-\infty, \lambda]$, for any $\Gamma$, where of course $\sigma$
has to be an irrational cohomology class. This has been shown to
be the case at least in the discretization of the
magnetic Schr\"odinger operator, called the Harper operator, and
when $\Gamma = \mathbb Z^2$, cf. \cite{BeSi}, \cite{CEY}, \cite{Last}.
There have been other interesting papers on gaps in the spectrum
of magnetic Schr\"odinger operators on covering spaces.
One such is \cite{KaPe} where the asymptotics of the size
of gaps in the spectrum are obtained on hyperbolic spaces,
as the mid-point of the gap tends to infinity, generalizing
known results on Euclidean spaces.
In two dimensions, magnetic Schr\"odinger operators on covering spaces
are the Hamiltonians in the model for the integer quantum Hall effect
when the covering space is the Euclidean plane \cite{Bel+E+S},
and for the fractional quantum Hall effect when the covering space
is the hyperbolic plane \cite{MM}, \cite{CHMM}.
The standard assumption made is that
the Fermi energy level lies in a gap of the spectrum of the
Hamiltonian (which however can be relaxed to the assumption that it lies
in a gap in extended states of the Hamiltonian after further analysis).

In this paper, we study the semiclassical approximation, as the
coupling  constant $\mu$
approaches zero, of the  Hamiltonians $H_{A,V}(\mu)$, for a
class of {\em Morse type} potentials $V$,
which include all functions $V=|df|^2$, where $f$ is a
$\Gamma$-invariant Morse function on
$\M$, that is
$f$ is the lift to $\M$ of a Morse function on $M$.  We show that the
spectrum  of this
operator, is approximated by the union of the spectra of  model
operators which are defined
near the critical points of the  Morse potential $V$. In particular,
we are able to deduce
the existence of an arbitrarily large number of gaps in the spectrum
of these Hamiltonians, for $\mu$
sufficiently small, a fact which is of crucial importance in the
study of the quantum Hall effect. We adapt the $L^2$ version of
semiclassical approximation of Witten \cite{Wi},
\cite{CFKS}, \cite{Sh} to our context.

The paper is organized as follows. We first give a summary of our
main results. Then we
recall some preliminary material on projective unitary
representations and the associated
von Neumann algebra of operators commuting with this algebra,
together with the von Neumann
trace and von Neumann dimension function.  In the next section, we
establish the first main
theorem on the existence of  spectral gaps for large values of the
coupling parameter or
equivalently for  small values of the coupling constant $\mu$. Here
we use the method  of
semiclassical approximation.

\subsection*{Summary of main results}

In this paper, we will assume that the potential $V$ is a smooth
$\Gamma$-invariant function which
satisfies the following {\em Morse type condition:}
$V(x) \ge 0$ for all $x\in \M$. Also if $V(x_0) = 0$ for some $x_0$
in $\M$, then
there is a  positive constant $c$ such that $V(x)\ge c|x-x_0|^2$ for all
$x$ in a neighborhood of $x_0$. We will also assume that $V$ has at
least one zero point.

We remark that all functions $V=|df|^2$, where $|df|$
denotes the pointwise norm of the differential of a $\Gamma$-invariant Morse
function $f$ on $\M$, are examples of Morse type potentials.

We next enunciate the principle on which this paper is based.
Associated to each
Hamiltonian $H_{A,V}(\mu)$, there is a {\em  model operator} $K$ (cf.
section 2)  which has a discrete spectrum.
It is defined as  a   direct sum of harmonic oscillators,
associated with the potential wells of $V$ in a fundamental domain
of $\Gamma$ in $\widetilde M$.

\noindent{\bf Semiclassical approximation principle:} {\em  Let $V$
be a Morse type
potential. Then in the semiclassical limit as the coupling constant
$\mu$ goes to
zero, the spectrum of the Hamiltonian $H_{A,V}(\mu)$ ``tears
up'' into bands which are located near the eigenvalues of the
associated model operator $K$.
}

The following main theorem establishes the existence of arbitrarily
large number of gaps in the spectrum of
$H_{A,V}(\mu)$ whenever $\mu$ is sufficiently small. The proof of
this theorem uses an
analogue of Witten's semiclassical approximation technique for
proving the Morse
inequalities. The
$L^2$-analogue (in the absence of a magnetic field) was proved by
Shubin \cite{Sh}, also
Burghelea et al. \cite{BFKM}. We modify the proof in \cite{Sh}
to obtain the result in section 2.
The physical explanation for the appearance of gaps in the spectrum
$H_{A,V}(\mu)$ is that the potential wells get deeper as $\mu\to 0$
and the atoms get (asymptotically) isolated, so that the energy
levels of $H_{A,V}(\mu)$ are approximated by those of the
corresponding model operator $K$.

\begin{thm*}[Existence of spectral gaps]
Let $V$ be a Morse type potential. If $E\in \C$ is such that
$E\not\in {\rm spec(K)}$, then
there exists $\mu_0 >0$ such that for all $\mu \in (0, \mu_0)$,
$E$ is in the resolvent set of $H_{A,V}(\mu)$. If in addition $E$ is
real and lies between two
eigenvalues of $K$, then $E$ is in
a spectral gap of $H_{A,V}(\mu)$.

Since the spacing between  the eigenvalues of $K$ is bounded below,
it follows that there
exists  arbitrarily large number of gaps in the spectrum of
$H_{A,V}(\mu)$ whenever the coupling constant $\mu$ is
sufficiently small.
\end{thm*}

\begin{rem*} We observe that the operator $H_A + \mu^{-2}V =
\mu^{-1}H_{A,V}(\mu)$ also
has arbitrarily large number of gaps in its spectrum whenever the
coupling constant $\mu$ is
sufficiently small. Analogous results in the special case of
Euclidean space were obtained
in \cite{Nak+Bel}.
\end{rem*}

\section{Preliminaries}

Let $M$ be a compact connected Riemannian manifold, $\Gamma$ be its
fundamental group and
$\widetilde M$ be its universal cover, i.e. one has the principal bundle
$ \Gamma\to \widetilde M\overset{p}{\to} M.$ To make the paper self-contained,
we include preliminary material, some of which may not be new, cf. 
\cite{Bel}, \cite{BrSu},
\cite{CHMM}, \cite{MM}, \cite{Ma}.

\subsection{\bfseries Projective action, or magnetic translations}

Let $\omega$ be a closed real-valued 2-form on $M$ such that $B =p^* \omega$ is
\emph{exact}. So $B=dA$ where $A$ is a 1-form on $\widetilde M$. We
will assume $A$
without loss of generality that $A$ is  real-valued too. Define
$\nabla=d+\,iA$.  Then
$\nabla$ is a Hermitian connection on the trivial line bundle  over
$\widetilde{M}$
with the curvature
$(\nabla)^2=i\, B$. The connection $\nabla$ defines  a projective
action of $\Gamma$
on $L^2$ functions as follows.

Observe that since $B$ is
$\Gamma$-invariant, one has
$ 0=\gamma^*B-B=d(\gamma^*A-A)\quad
\forall\gamma \in\Gamma$. So $\gamma^*A-A$ is a closed 1-form on the
simply connected
manifold $\widetilde{M}$, therefore
\[
\gamma^*A-A=d\psi_\gamma,\quad\forall\gamma\in\Gamma,
\]
where $\psi_\gamma$ is a smooth
function on $\widetilde{M}$. It is defined up to an additive constant,
so we can assume in addition that  it satisfies
the following normalization condition:
\begin{itemize}
\item $\psi_\gamma(x_0)=0$ for a fixed $x_0\in\widetilde{M},\quad
\forall\gamma\in\Gamma$.
\end{itemize}
It follows that $\psi_\gamma$ is real-valued and $\psi_e(x)\equiv 0$, where
$e$ denotes the neutral element of $\Gamma$.
It is also easy to check that

\begin{itemize}
\item $\psi_\gamma(x)+\psi_{\gamma'}(\gamma
x)-\psi_{\gamma'\gamma}(x)$  is independent of
$x\in\widetilde{M},\quad \forall \gamma,\gamma'\in\Gamma$.
     \end{itemize}

Then $\sigma(\gamma,\gamma')=\exp(-i\psi_\gamma(\gamma'\cdot x_0))$  defines a
{\em multiplier} on $\Gamma$ i.e. $\sigma:\Gamma\times\Gamma\to U(1)$
satisfies
\begin{itemize}
\item ${\sigma}(\gamma,e) = {\sigma}(e,\gamma)=1,\quad\forall\
\gamma\in\Gamma$;

\item ${\sigma}(\gamma_1,\gamma_2)
{\sigma}(\gamma_1\gamma_2, \gamma_3)=
{\sigma}(\gamma_1,\gamma_2\gamma_3)
{\sigma}(\gamma_2,\gamma_3),\quad \forall \gamma_1, \gamma_2,
\gamma_3\in \Gamma$
\quad ({\em the cocycle relation}).
\end{itemize}

It follows from these relations that
$\sigma(\gamma,\gamma^{-1})=\sigma(\gamma^{-1},\gamma)$.

The complex conjugate multiplier
$\bar\sigma(\gamma,\gamma')=\exp(i\psi_\gamma(\gamma'\cdot x_0))$
also satisfies the same relations.

For $u\in L^2(\widetilde M)$ and $\gamma\in\Gamma$ define
$$
\quad U_\gamma u = (\gamma^{-1})^* u, \quad S_\gamma u =
\exp(-i\psi_\gamma)\,u.
$$
Then the operators $T_\gamma=U_\gamma\circ S_\gamma$ satisfy
$$
T_e={\rm Id},
     \quad T_{\gamma_1} T_{\gamma_2}
= {\sigma}(\gamma_1,\gamma_2) T_{\gamma_1 \gamma_2},
$$
for all $\gamma_1, \gamma_2 \in \Gamma$. In this case one says that the map
$T : \Gamma\to {\U} (L^2(\M))$, $\gamma\mapsto T_\gamma$, is a projective
$(\Gamma, {\sigma})$-unitary representation, where for any Hilbert space $\H$
     we denote by
$\U(\H)$ the group of all unitary operators in $\H$. In other words
one says that the map $\gamma\mapsto T_\gamma$ defines a
$(\Gamma,\sigma)$-{\em action}
in~$\H$.

It is also easy to check that the adjoint operator to $T_\gamma$ in
$L^2(\widetilde M)$
(with respect to a smooth $\Gamma$-invariant measure) is
$$
T_\gamma^*=\bar\sigma(\gamma,\gamma^{-1})T_{\gamma^{-1}}.
$$

The operators $T_\gamma$ are also called
{\em magnetic translations}.

\subsection{Twisted group algebras}
Denote by $\ell^2(\Gamma)$ the standard Hilbert space of
complex-valued $L^2$-functions
on the discrete group $\Gamma$.
We will use a left $(\Gamma, \bar\sigma)$-action on
$\ell^2(\Gamma)$ (or, equivalently, a $(\Gamma, \bar\sigma)$-unitary
representation
in $\ell^2(\Gamma)$) which is given explicitly by
$$
T_\gamma^L f(\gamma') = f(\gamma^{-1}\gamma')
\bar\sigma(\gamma, \gamma^{-1}\gamma'), \qquad \gamma, \gamma' \in \Gamma.
$$
It is easy to see that this is indeed a $(\Gamma,\bar\sigma)$-action, i.e.
$$
T_e^L={\rm Id} \quad {\rm and} \quad
T_{\gamma_1}^L T_{\gamma_2}^L=
\bar\sigma(\gamma_1,\gamma_2)T_{\gamma_1 \gamma_2}^L, \quad \forall
\gamma_1,\gamma_2\in\Gamma.
$$
Also
$$
(T_\gamma^L)^*=\sigma(\gamma, \gamma^{-1})T_{\gamma^{-1}}^L.
$$

Let
$$
{\mathcal A}^R(\Gamma, \sigma) =
\Big\{ A \in \B(\ell^2(\Gamma)): [T_\gamma^L, A] = 0, \quad \forall \gamma \in
\Gamma\Big\}
$$
be the commutant of the left $(\Gamma, \bar\sigma)$-action on
$\ell^2(\Gamma)$.  Here by $\B(\H)$ we denote the algebra of all
bounded linear operators
in a Hilbert space $\H$.
By the general theory, ${\mathcal A}^R(\Gamma, \sigma)$ is a von
Neumann algebra
and is known as the {\em (right) twisted group von Neumann algebra}.
It can also be realized as follows. Let us define the  following operators
in $\ell^2(\Gamma)$:
$$
T_\gamma^R f(\gamma') = f(\gamma'\gamma) \sigma(\gamma', \gamma),
\qquad  \gamma, \gamma' \in \Gamma.
$$
It is easy to check that they form a right $(\Gamma,\sigma)$-action
in $\ell^2(\Gamma)$ i.e.
$$
T_e^R={\rm Id} \quad {\rm and} \quad
T_{\gamma_1}^R  T_{\gamma_2}^R=
\sigma(\gamma_1,\gamma_2)T_{\gamma_1 \gamma_2}^R, \quad \forall
\gamma_1,\gamma_2\in\Gamma,
$$
and also
$$
(T_{\gamma}^R)^*=\bar\sigma(\gamma,\gamma^{-1}) T_{\gamma^{-1}}^R.
$$
This action commutes with the left $(\Gamma,\bar\sigma)$-action
defined above i.e.
$$
T_\gamma^L  T_{\gamma'}^R= T_{\gamma'}^R T_\gamma^L,
\quad \forall \gamma,\gamma'\in\Gamma.
$$
It can be shown that  the von Neumann algebra $\A^R(\Gamma,\sigma)$
is generated by the
operators
$\{T_{\gamma}^R\}_{\gamma\in\Gamma}$
(see e.g.\ a similar argument in \cite{Shubin}).

Similarly we can introduce  a von Neumann algebra
$$
{\mathcal A}^L(\Gamma, \bar\sigma) =
\Big\{ A \in \B(\ell^2(\Gamma)): [T_\gamma^R, A] = 0, \quad \forall \gamma \in
\Gamma\Big\}.
$$
We will refer to it as {\em (left) twisted group von Neumann algebra}.
It is generated by the operators
$\{T_{\gamma}^L\}_{\gamma\in\Gamma}$, and it is the commutant
of ${\mathcal A}^R(\Gamma, \sigma)$.

Let us define a twisted group algebra $\C(\Gamma,\sigma)$ which consists of
complex valued functions with finite support on $\Gamma$ and with the
twisted convolution
operation
$$
(f*g)(\gamma)=\sum_{\gamma_1,\gamma_2:\gamma_1\gamma_2=
\gamma}f(\gamma_1)g(\gamma_2)\sigma(\gamma_1,\gamma_2).
$$
The basis of $\C(\Gamma,\sigma)$ as a vector space is formed by
$\delta$-functions
$\{\delta_\gamma\}_{\gamma\in\Gamma}$, $\delta_\gamma(\gamma')=1$ if
$\gamma=\gamma'$
and $0$ otherwise. We have
$$
\delta_{\gamma_1} *
\delta_{\gamma_2}=\sigma(\gamma_1,\gamma_2)\delta_{\gamma_1\gamma_2}.
$$
Associativity of this multiplication is equivalent to the cocycle condition.

Note also that the $\delta$-functions $\{\delta_\gamma\}_{\gamma\in\Gamma}$
form an orthonormal basis in $\ell^2(\Gamma)$. It is easy to check that
$$
T_\gamma^L\delta_{\gamma'}=\bar\sigma(\gamma,\gamma')\delta_{\gamma\gamma'},
\quad
T_\gamma^R\delta_{\gamma'}=
\sigma(\gamma'\gamma^{-1},\gamma)\delta_{\gamma'\gamma{-1}}.
$$

It is clear that the correspondences $\delta_\gamma\mapsto T^L_{\gamma}$
and $\delta_\gamma\mapsto  T^R_{\gamma}$ define representations
of $\C(\Gamma,\bar\sigma)$ and $\C(\Gamma,\sigma)$ respectively. In both cases
the weak closure of the image of the twisted group algebra coincides
with the corresponding
von Neumann algebra ($\A^L(\Gamma,\bar\sigma)$ and
$\A^R(\Gamma,\sigma)$ respectively). The
corresponding norm closures are so called {\em reduced twisted group}
$C^*$-{\em algebras} which are denoted $C^*_r(\Gamma,\bar\sigma)$ and
$C^*_r(\Gamma,\sigma)$
respectively.

The von Neumann algebras $\A^L(\Gamma,\bar\sigma)$ and
$\A^R(\Gamma,\sigma)$ can be described
in terms of the matrix elements. For any $A\in\B(\ell^2(\Gamma))$ denote
$A_{\alpha,\beta}=(A\delta_\beta,\delta_\alpha)$ (which is a matrix
element of $A$).
Then repeating standard arguments (given in a similar situation e.g.\
in \cite{Shubin}) we
can prove that for any $A\in\B(\ell^2(\Gamma))$ the inclusion $A\in
\A^R(\Gamma,\sigma)$
is equivalent to the relations
$$
A_{\gamma x,\gamma y}=\bar\sigma(\gamma,x)\sigma(\gamma,y)A_{x,y}\;,\quad
\forall x,y,\gamma\in\Gamma.
$$
In particular, we have for any $A\in \A^R(\Gamma,\sigma)$
$$
A_{\gamma x,\gamma x}=A_{x,x}\;, \quad \forall x,\gamma\in\Gamma.
$$
Similarly, for any $A\in\B(\ell^2(\Gamma))$ the inclusion $A\in
\A^L(\Gamma,\bar\sigma)$
is equivalent to the relations
$$
A_{x\gamma,y\gamma}=\bar\sigma(x,\gamma)\sigma(y,\gamma)A_{x,y}\;,\quad
\forall x,y,\gamma\in\Gamma.
$$
In particular, we have
$$
A_{x\gamma ,x\gamma }=A_{x,x}\;, \quad \forall x,\gamma\in\Gamma,
$$
for any $A\in \A^L(\Gamma,\bar\sigma)$.

A finite {\em von Neumann trace}
$\tr_{\Gamma,\bar\sigma}:\A^L(\Gamma,\bar\sigma)\to\C$
is defined by the formula
$$
\tr_{\Gamma,\bar\sigma} A=(A\delta_e,\delta_e).
$$
We can also write $\tr_{\Gamma,\bar\sigma} A=A_{\gamma,\gamma}=
\left(A\delta_\gamma,\delta_\gamma\right)$  for any $\gamma\in
\Gamma$ because the right hand
side does not depend of
$\gamma$.

A finite von Neumann trace $\tr_{\Gamma,\sigma}:\A^R(\Gamma,\sigma)\to\C$
is defined by the same formula, so we will denote by $\tr_\Gamma$ any
of these traces.

Let $\mathcal H$ denote an infinite dimensional complex Hilbert space.
Then the Hilbert tensor product $\ell^2(\Gamma)\otimes \mathcal H$ is both
$(\Gamma, \bar\sigma)$-module and $(\Gamma, \sigma)$-module under the actions
$\gamma \mapsto T_\gamma^L \otimes 1$ and $\gamma \mapsto  T_\gamma^R
\otimes 1$
respectively. Let ${\mathcal A}^L_{\mathcal H}(\Gamma,\bar\sigma)$
and ${\mathcal A}^R_{\mathcal H}(\Gamma,\sigma)$ denote the von Neumann
algebras in $\ell^2(\Gamma)\otimes \mathcal H$ which are
commutants
of the $(\Gamma, \sigma)$- and $(\Gamma, \bar\sigma)$-actions respectively.
Clearly
${\mathcal A}^L_{\mathcal H}(\Gamma,\bar\sigma)
\cong {\mathcal A}^L(\Gamma,\bar\sigma) \otimes \B(\mathcal H)$
and
${\mathcal A}^R_{\mathcal H}(\Gamma,\sigma)
\cong {\mathcal A}^R(\Gamma,\sigma) \otimes \B(\mathcal H)$
in the usual sense of von Neumann algebra tensor products.

Define the semifinite
tensor product trace
$\Tr_{\Gamma} = \tr_{\Gamma} \otimes \Tr$
on each of the algebras ${\mathcal A}^L_{\mathcal H}(\Gamma,\bar\sigma)$
and ${\mathcal A}^R_{\mathcal H}(\Gamma,\sigma)$.
Here $\Tr$ denotes the standard (semi-finite) trace on $\B(\H)$.

Let us recall that a closed linear subspace $V\subset \widetilde \H$ is
{\it affiliated} to a von Neumann algebra $\A$ of operators in a Hilbert space
$\widetilde \H$ if
$P_V\in\A$ where $P_V$ is the orthogonal projection on $V$ in $\widetilde \H$.
This is equivalent to saying that $V$ is invariant under the
commutant $\A'$ of $\A$
in $\widetilde \H$.

The {\em von Neumann dimension} $\dim_{\Gamma}$ of a
closed subspace $V$ of $\ell^2(\Gamma)\otimes \mathcal H$ that is
invariant under
$\{T_\gamma^L\otimes 1,\;\gamma\in\Gamma\}$ (or, equivalently, affiliated to
$\A_\H^R(\Gamma,\sigma)$) is defined as
$$
\dim_{\Gamma}(V) = \Tr_{\Gamma}(P_V).
$$
The same formula is used for the subspaces which are
invariant under
$\{T_\gamma^R\otimes 1,\;\gamma\in\Gamma\}$ (or, equivalently, affiliated to
$\A_\H^L(\Gamma,\bar\sigma)$).

Also the {\em von Neumann rank} of an operator
$Q\in{\mathcal A}^R_{\mathcal H}(\Gamma,\sigma)$ or
$Q\in{\mathcal A}^L_{\mathcal H}(\Gamma,\bar\sigma)$ is defined as
$$
\rank_\Gamma(Q) = \dim_{\Gamma}(\overline{{\rm Range}(Q)}),
$$
where the bar over ${\rm Range}(Q)$ means closure.

\subsection{\bfseries Magnetic Hamiltonians} The {\em magnetic
Laplacian} on $L^2(\M)$ is
defined as
\[
H_A = \nabla^*\nabla = \left(d+iA\right)^*
\left(d+iA\right)
\]
and more generally, the {\em magnetic Hamiltonian} or
{\em Magnetic Schr\"odinger operator} is defined as
\[
H_{A,V}(\mu) = \mu H_A + \mu^{-1}V,
\]
where $V$ is any $\Gamma$-invariant smooth function on $\M$.
The Hamiltonian $H=H_{A,V}(\mu)$ is a self adjoint second order
elliptic differential
operator. It commutes with the  magnetic  translations $T_\gamma$ (for all
$\gamma\in\Gamma$), i.e.  with the $(\Gamma,{\sigma})$-action which
was defined above. To see this note first that the operators
$U_\gamma=(\gamma^{-1})^*$ and $S_\gamma$ (the multiplication by
$\exp(-i\psi_\gamma)$)
are defined not only on scalar functions but also on 1-forms
(and actually on $p$-forms for any $p\ge 0$) on $\widetilde M$.
Hence the magnetic translations $T_\gamma$ are well defined on forms as well.
The operators $T_\gamma$ are obviously unitary on the $L^2$ spaces of forms,
where the $L^2$ structure is defined by the fixed $\Gamma$-invariant metric on
$\widetilde M$.  An easy calculation shows that $T_\gamma
\nabla=\nabla T_\gamma$
on scalar functions. By taking adjoint operators we obtain
$T_\gamma \nabla^*=\nabla^* T_\gamma$ on 1-forms. Therefore $T_\gamma
H_A=H_A T_\gamma$
on functions. Since obviously $T_\gamma V=V T_\gamma$, we see that
$H_{A,V}(\mu)$
commutes with $T_\gamma$ for all $\gamma$.

In dimension 2, it is the spectrum and the spectral projections of the magnetic
Schr\"odinger operator that are of
fundamental importance to the study of the quantum Hall effect. We
remark that it is
virtually impossible
to explicitly compute the spectrum of $H_{A,V}$ for arbitrary $V$ which is
$\Gamma$-invariant, even in 2 dimensions and for simplest manifolds.
Nevertheless,
Comtet and Houston
\cite{Comtet+H} computed the
spectrum of $H_A$ (with $V=0$) on the hyperbolic plane  with the magnetic
potential $A = \theta y^{-1}dx$ (which corresponds to the constant magnetic
field), where we can  assume without loss of generality that $\theta>0$.
This spectrum is  the
union of a finite number of eigenvalues
$
\left\{(2k+1)\,\theta-k(k+1) : k=0,1,2,\dots <\theta-\frac{1}{2}\right\}
$ and the continuous spectrum
$
\left[\frac{1}{4} + \theta^2, \infty \right)$.

Since $H$ commutes with the $(\Gamma,{\sigma})$-action, it follows by
the spectral
mapping theorem that  the spectral projections of
$H$, $E_\lambda = \chi_{(-\infty,\lambda]}(H)$ are
bounded operators on
$L^2(\M)$ that also commute with the $(\Gamma,{\sigma})$-action
i.e. $T_\gamma E_\lambda= E_\lambda T_\gamma,\quad\forall\ \gamma\in\Gamma$.
The {\em commutant} of the $(\Gamma, {\sigma})$-action is a von Neumann algebra
\[
{\mathcal U}_{\M}(\Gamma,\bar\sigma) = \left\{Q\in \B(L^2(\M)) :
T_\gamma Q = Q T_\gamma,\quad\forall\ \gamma\in\Gamma\right\}.
\]
To characterize the Schwartz kernels $k_Q(x,y)$ of the operators
$Q\in {\mathcal U}_{\M}(\Gamma,\bar\sigma)$ note that the relation
$T_\gamma Q=QT_\gamma$
can be rewritten in the form
$$
e^{i\psi_\gamma(x)} k_Q(\gamma x,\gamma y)
e^{-i\psi_\gamma(y)} = k_Q(x,y), \quad
\forall x,y\in\M\quad \forall \gamma\in\Gamma,
$$
so $Q\in {\mathcal U}_{\M}(\Gamma,\bar\sigma)$ if and only if this holds
for all $\gamma\in\Gamma$. In particular, in this case
$k_Q(x,x)$ is $\Gamma$-invariant. For the spectral projections of $H$
we also have
$E_\lambda\in {\mathcal U}_{\M}(\Gamma,\bar\sigma)$, so the
corresponding Schwartz kernels also
satisfy the relations above. Note that the Schwartz kernels of
$E_\lambda$ are in
$C^\infty(\widetilde M\times\widetilde M)$.

To define a natural trace on ${\mathcal U}_{\M}(\Gamma,\bar\sigma)$
we will construct
an isomorphism of this algebra with the von Neumann algebra
$\A^L_\H(\Gamma,\bar\sigma)$,
where $\H=L^2(\F)$.

Let $\mathcal{F}$ be a fundamental domain for the $\Gamma$-action on
$\M$. By choosing a connected fundamental domain $\mathcal{F}$ for
the action of
$\Gamma$ on $\M$, we can
define a $(\Gamma, \sigma)$-equivariant isometry
$$
{\bf U} : L^2(\M)\cong \ell^2(\Gamma)\otimes L^2(\mathcal{F}) \leqno (1.1)
$$
as follows. Let $i: \mathcal{F} \to \M$ denote the inclusion map.
Define
$$
{\bf U} (\phi) = \sum_{\gamma\in \Gamma}\delta_\gamma\otimes
i^*(T_\gamma \phi),
\qquad  \phi \in L^2(\M).
$$
\begin{lemma} The map
$\; {\bf U} :  L^2(\M)\to  \ell^2(\Gamma)\otimes L^2(\mathcal{F})\;$ defined
above is a
     $(\Gamma, \sigma)$-equivariant isometry, where the
$(\Gamma, \sigma)$-action is given by the operators
$T_\gamma$ and $T^R_\gamma\otimes 1$ on the spaces $L^2(\M)$ and
  $\ell^2(\Gamma)\otimes L^2(\mathcal{F})$ respectively.
\end{lemma}

\begin{proof} Given $\phi \in L^2(\M)$, we compute

\begin{align*}
{\bf U}(T_\gamma\phi)
&=\sum_{\gamma'\in\Gamma}\delta_{\gamma'}\otimes i^*(T_{\gamma'}T_\gamma\phi)
=\sum_{\gamma'\in\Gamma}\sigma(\gamma',\gamma)\delta_{\gamma'}\otimes
i^*(T_{\gamma'\gamma}\phi)\\
&=\sum_{\gamma'\in\Gamma}\sigma(\gamma'\gamma^{-1},\gamma)\delta_{\gamma'\gamma^{-1}}
\otimes i^*(T_{\gamma'}\phi)
=(T^R_\gamma\otimes 1){\bf U}\phi,
\end{align*}
which proves that $\bf U$ is a $(\Gamma,\sigma)$-equivariant map. It
is straightforward
to check that that $\bf U$ is an isometry.
\end{proof}

Since $\mathcal{U}_{\M}(\Gamma,\bar\sigma)$ is the commutant of
$\{T_\gamma\}_{\gamma\in\Gamma}$,
and $\A^L_\H(\Gamma,\bar\sigma)$, with $\H=L^2(\F)$, is the commutant of
$\{T^R_\gamma\otimes 1\}_{\gamma\in\Gamma}$, we see that $\bf U$
induces an isomorphism
of von Neumann algebras $\mathcal{U}_{\M}(\Gamma,\bar\sigma)$ and
$\A^L_\H(\Gamma,\bar\sigma)$.
Therefore we can transfer the trace $\Tr_\Gamma$ from
$\A^L_\H(\Gamma,\bar\sigma)$ to
$\mathcal{U}_{\M}(\Gamma,\bar\sigma)$. The result will be a
semifinite $\Gamma$-trace on
$\mathcal{U}_{\M}(\Gamma,\bar\sigma)$ which we will still denote $\Tr_\Gamma$.

It is easy to check that for any $Q\in
\mathcal{U}_{\M}(\Gamma,\bar\sigma)$ with
a finite $\Gamma$-trace and a continuous Schwartz kernel $k_Q$ we have
$$
\Tr_\Gamma Q=\int_\F k_Q(x,x)dx
$$
where $dx$ means the $\Gamma$-invariant measure. An important
particular case is
a spectral projection $E_\lambda$ of the magnetic Schr\"odinger
operator $H=H_{A,V}$
as considered above. The projection $E_\lambda$  has a finite $\Gamma$-trace
and a $C^\infty$ Schwartz kernel. Therefore we can define a {\em
spectral density function}
$$
N_\Gamma(\lambda;H)=\Tr_\Gamma E_\lambda,
$$
which is finite for all $\lambda\in \R$. It is easy to see that
$\lambda\mapsto N_\Gamma(\lambda;H)$
is a non-decreasing function, and the
spectrum of
$H$ can be reconstructed  as the set of its points of growth, i.e.
$$
\spec (H)=\{\lambda\in\R:
N_\Gamma(\lambda+\eps;H)-N_\Gamma(\lambda-\eps;H)>0, \ \forall
\eps>0\}.
$$

The {\em von Neumann dimension} $\dim_\Gamma$ of a
closed subspace $V$ of $L^2(\widetilde M)$ is well defined if $V$ is
invariant under
$T_\gamma,\;\forall \gamma\in\Gamma$. Also the {\em von Neumann rank}
of an operator
$Q \in {\mathcal U}_{\M}(\Gamma,\bar\sigma)$ is well defined.

\subsection{Gauge invariance}

If we make another choice of vector potential $A'$ such that $dA' = B$, then
it follows that $A'-A$ is a closed $1$-form on a simply connected manifold
$\M$, and therefore it is exact, i.e. $A'=A+ d\phi$, where $\phi\in
C^\infty(\widetilde M)$.
We will always assume that $\phi$ is normalized  by the condition
$\phi(x_0)=0$. It follows that
$\phi$ is real-valued.

It then
follows that  the connection $\nabla = d+iA$ gets unitarily
conjugated into a new connection
$\nabla' = d+iA' = e^{-i\phi}\nabla   e^{i\phi}$.  Therefore
${\nabla'}^*=e^{-i\phi}\nabla^*  e^{i\phi}$ and $H'=e^{-i\phi}H
e^{i\phi}$, where
$H=H_{A,V}$, $H'=H_{A',V}$,and $V$ is $\Gamma$-invariant.
In particular, $H'$ and $H$ are unitarily equivalent.
Let us repeat the constructions of the previous subsections with $A$
replaced by $A'$
indicating relations of the modified objects with the old ones.

Define the function  $\psi'_\gamma$ from
$d\psi'_\gamma=\gamma^*A'-A',$ $\psi'_\gamma(x_0)=0$
(with the same point $x_0\in M$ as above). Then
$$
\psi'_\gamma=\psi_\gamma+\gamma^*\phi-\phi-\phi(\gamma x_0).
$$
The new cocycle will be
$$
\sigma'(\gamma_1,\gamma_2)=\exp\left(-i\psi'_{\gamma_1}(\gamma_2
x_0)\right) =\sigma(\gamma_1,\gamma_2)
\exp\left(-i[\phi(\gamma_1\gamma_2 x_0)-\phi(\gamma_1
x_0)-\phi(\gamma_2 x_0)]\right).
$$
The modified magnetic translations are defined by
$T'_\gamma=U'_\gamma S'_\gamma$, where
$U'_\gamma=U_\gamma=(\gamma^{-1})^*$ and
$S'_\gamma=\exp(-i\psi'_\gamma)$. Then
$T'_{\gamma_1}T'_{\gamma_2}=\sigma'(\gamma_1,\gamma_2)T'_{\gamma_1 \gamma_2}$.

Clearly $H'$ commutes with the modified magnetic translations
$T'_\gamma,\ \forall\gamma\in\Gamma$.
The relation between old and new magnetic translations is
$$
T'_\gamma=e^{i\phi(\gamma\cdot x_0)}\left(e^{-i\phi}T_\gamma e^{i\phi}\right),
$$
which is again the same unitary conjugation up to a constant unitary factor.

Now we can introduce a von Neumann algebra
$$
\U_{\widetilde M}(\Gamma,\bar\sigma')=\{Q'\in\B(L^2({\widetilde M})):
T'_\gamma Q'=Q' T'_\gamma,\
\forall \gamma\in\Gamma\}.
$$
Clearly the map
$$
\U_{\widetilde M}(\Gamma,\bar\sigma)\longrightarrow \U_{\widetilde
M}(\Gamma,\bar\sigma'),\qquad
Q\longmapsto Q'=e^{-i\phi}Qe^{i\phi},
$$
is an isometric $\star$-isomorphism of von Neumann algebras. It is
easy to see that this isomorphism
preserves the $\Gamma$-trace which is defined on both algebras. If an operator
$Q$ has a smooth Schwartz kernel $k_Q(x,y)$ and a finite $\Gamma$-trace, then
$$
k_{Q'}(x,y)=e^{-i\phi(x)}k_Q(x,y)e^{i\phi(y)}, \quad\forall x,y\in
\widetilde M,
$$
hence
$$
k_{Q'}(x,x)=k_Q(x,x), \quad \forall x\in \widetilde M,
$$
so the equality $\Tr_\Gamma Q'=\Tr_\Gamma Q$ follows from the
expression of the $\Gamma$-traces
in terms of kernels.

If $E_\lambda$ and $E'_\lambda$
are spectral projections of $H$ and $H'=e^{-i\phi}He^{i\phi}$
respectively, then clearly
$E'_\lambda=e^{-i\phi}E_\lambda e^{i\phi}$, therefore $\Tr_\Gamma
E'_\lambda=\Tr_\Gamma E_\lambda$.
This means that the spectral density functions of $H$ and $H'$
coincide, i.e. the spectral density
function is gauge invariant.

\section{Semiclassical approximation and the existence of spectral gaps}
We will
study an $L^2$-version of semiclassical approximation,  which is
similar to the ones
which appear when we take the Witten deformation of the de Rham
complex and consider the
corresponding Laplacian (cf. \cite{{Wi},{CFKS},{HS}}. For the case of
the algebra
corresponding to the regular representation of $\pi_1(M)$, such
asymptotics were first
proved in \cite{Sh}, see also \cite{BFKM} for a related semiclassical
approximation technique. The proofs given in \cite{Sh} will be adapted
to work in the more general situation that we need in this section.

Recall that $H = H_{A, V}(\mu)=\mu H_A+\mu^{-1}V$ is a second order
differential
operator acting on $L^2(\widetilde M)$, such that it commutes 
with
the projective unitary
$(\Gamma, \sigma)$-action on 
$L^2(\widetilde M)$, given by the
magnetic translations
$\{T_\gamma,\ 
\gamma\in\Gamma\}$.  Moreover, note that
$H$ is a  second order 
elliptic operator with a positive principal symbol,
order operator,
$V$ is a non-negative potential function on $M$ which 
has only
nondegenerate zeroes and
$\mu>0$ is a small 
parameter.

Actually the results will not change if we add to $H$ 
any
$\Gamma$-invariant zeroth order operator,
i.e. multiplication by 
a smooth $\Gamma$-invariant function.

Let us recall that 
$N_{\Gamma} (\lambda;\,H)$ denote the von
Neumann spectral density 
function  of the
operator $H$ which can be defined 
as
$$
N_{\Gamma}(\lambda, H) = \Tr_{\Gamma} 
E_\lambda(H)=
\Tr_{\Gamma}\left(\chi_{(-\infty,\lambda]}(H)\right),
$$
where 
$\chi_F$ means the characteristic function of a subset $F\subset 
\R$.

Let us choose a fundamental domain $\F\subset\widetilde M$ so 
that
there is no zeros of $V$
on the boundary of $\F$. This is equivalent to saying that  the translations
$\{\gamma\F,\;\gamma\in\Gamma\}$ cover  the set $V^{-1}(0)$ (the set
of all zeros of $V$).
Let $V^{-1}(0)\cap\F=\{\bar x_j|\,j=1,\dots,N\}$ be the set of all
zeros of $V$ in $\F$;
$\bar x_i\ne\bar x_j$ if $i\ne j$.
     Let $K$ denote the {\em model operator} of $H$ (cf.
\cite{Sh}), which is obtained as a direct sum of quadratic parts of $H$ in all
points $\bar x_1,\dots,\bar x_N$.
More precisely,
$$
K = \oplus_{1\le j\le N} K_j,
$$
where
$K_j$ is an unbounded self-adjoint operator in $L^2({\mathbb R}^n)$ which
corresponds to the zero $\bar x_j$. It is a quantum harmonic oscillator and has
a discrete spectrum.
We assume that we have fixed local coordinates
on $\M$ in a small neighborhood $B(\bar x_j, r)$ of
     $\bar x_j$ for every $j=1,\dots,N$. Then $K_j$ has the form
$$
K_{j}= H_{j}^{(2)}+ V^{(2)}_{j},
$$
where all the components are obtained from $H$ as follows. In the
fixed local coordinates on
$\M$ near $\bar x_j$, the second order term
$H_{j}^{(2)}$ is a homogeneous second order  differential operator
with constant
coefficients (without lower order terms) obtained by isolating the
second order terms in
the operator $H$ and freezing the coefficients of this operator at
$\bar x_j$. (Note that $H_j^{(2)}$ does not depend of $A$.)
The zeroth order term $V^{(2)}_{j}$
is  obtained by taking the quadratic part of $V$ in the chosen
coordinates near $\bar x_j$.

More explicitly,
$$
H_j^{(2)}=\sum_{i,k=1}^n g^{ik}(\bar x_j)\frac{\partial^2}{\partial 
x_i\partial x_k},
\qquad
V_j^{(2)}=
\frac{1}{2}\sum_{i,k=1}^n\frac{\partial^2 V}{\partial x_i\partial 
x_k}(\bar x_j)x_ix_k,
$$
where $(g^{ik})$ is the inverse matrix to the matrix of the 
Riemannian tensor $(g_{ik})$.

We will say that $H$ is {\it flat} near $\bar x_j$ if $H=K_j$ near $\bar x_j$.
(In particular, in this case we should have $A=0$ near $\bar x_j$.)

We will also need the operator
$$
K(\mu) = \oplus_{1\le j\le N} 
K_j(\mu),
$$
where
$$
K_j(\mu)=\mu H_{j}^{(2)}+\mu^{-1} V^{(2)}_{j}, 
\qquad \mu>0.
$$
It is easy to see that $K(\mu)$ has the same 
spectrum as $K=K(1)$.


Let $\{\alpha_p: p \in \mathbb N \}$ denote the set of all
eigenvalues of the model operator
$K$, $\alpha_p\not=\alpha_q$ for $p\not=q$, and $r_p$ denote the
multiplicity of
$\alpha_p$, i.e. $r_p=\dim \Ker(K-\alpha_pI)$ where $K$ is considered
in $(L^2(\R^n))^N$.
Denote by $N(\lambda;K)$ the distribution function of the eigenvalues
of $K$, i.e.
$N(\lambda;K)$ is the number of eigenvalues which are $\le\lambda$
(multiplicities counted).

The following is the main result of this paper.

\begin{thm}[Semiclassical Approximation]\label{T:semi}  For any $R>0$
there exist constants
$C>0$  and $\mu_0>0$ such that for any $\mu\in (0,\mu_0)$
and $\lambda\le R$, one has
$$
N(\lambda-C\mu^{1/5};K)\le N_\Gamma(\lambda;H)\le N(\lambda+C\mu^{1/5};K).
$$
Therefore
$${\rm spec}(H)\cap(-\infty,R]\subset\bigcup_{p=1}^\infty
(\alpha_p-C\mu^{1/5},\alpha_p+C\mu^{1/5})\,,$$
where $\{\alpha_p: p \in \mathbb N \}$ denotes the set of all
eigenvalues of the model operator
$K$.
Moreover for any $p=1,2,3,\dots$ {\it with}
$\alpha_p\in [-R,R]$ and any $\mu\in (0,\mu_0)$ one has
$$
N_{\Gamma}(\alpha_p+C\mu^{1/5};\,H) -
N_{\Gamma}(\alpha_p-C\mu^{1/5};\,H) =r_p
= N(\alpha_p+0;\,K)-N(\alpha_p-0;\,K).
$$
\end{thm}

This means that for small values of $\mu$, the spectrum of $H$
concentrates near the
eigenvalues of the model operator $K$, and for every such eigenvalue,
the von Neumann
dimension of the spectral subspace of the operator $H$, corresponding
to the part of the
spectra near the eigenvalue, is exactly equal to the usual multiplicity of this
eigenvalue of $K$.

\begin{proof}[Proof of Theorem on Existence of Spectral Gaps]
Let us 
assume that $\lambda \in \R$ and $\lambda \not\in {\rm spec}(K)$,
and 
then choose $R>\lambda$. Let
    $C>0$ and  $\mu_0>0$ be as in 
Theorem \ref{T:semi}. Then by taking
even smaller $\mu_0$ we will 
get
$$
    \lambda\not\in 
\bigcup_{p=1}^\infty
(\alpha_p-C\mu^{1/5},\alpha_p+C\mu^{1/5}), \quad 
\forall \mu\in(0,\mu_0).
$$
Therefore
$ \lambda\not\in {\rm 
spec}(H)$.
\end{proof}

The proof of Theorem 2.1 will be divided into 
2 parts: estimating
$N_\Gamma(\lambda,H)$
from below and from 
above.

\subsection{Estimate from below}

We will start by proving 
an estimate from below for
$N_\Gamma(\lambda;H)$ where 
$H=H(\mu)$,
$\lambda\le R$ with an arbitrarily fixed $R>0$, and 
$\mu\downarrow 0$.

Our main tool will be the standard variational 
principle for the
spectral density function.

\begin{lemma} 
[Variational principle] \label{L:var}
For every 
$\lambda\in\R$
\begin{equation}\label{E:var}
N_\Gamma(\lambda;H)=\sup\{\dim_\Gamma 
L\,|\; L\subset\Dom(H);\quad(Hf,f)\le
\lambda(f,f),\quad \forall f\in 
L\}.
\end{equation}
It is understood here that $L$ should be a 
closed
$(\Gamma,\sigma)$-invariant subspace in
$L^2(\widetilde M)$, 
i.e. closed subspace which is invariant under
all magnetic 
translations
$T_\gamma$, $\gamma\in\Gamma$.
\end{lemma}
A similar 
variational principle for the usual action of $\Gamma$
was used e.g. 
in \cite{ES,Sh}.

We will now describe an appropriate construction of 
a test space $L$.

Fix a function $J\in C_0^\infty(\R^n)$ such that 
$0\leq J\leq 1$,
$J(x)=1$ if $|x|\leq 1$, $J(x)=0$ if $|x|\geq 2$, 
and
$(1-J^2)^{1/2}\in C^\infty(\R^n)$.
Let us fix a
number
$\ka,\ 
0<\ka<1/2,$ which we shall choose later. For any $\mu>0$ 
define
$J^{(\mu)}(x)=J(\mu^{-\ka}x)$. This will be our standard 
cut-off function.
Let $J_j=J^{(\mu)}$ in the fixed coordinates near 
$\bar x_j$. Denote also
$J_{j,\gamma}=(\gamma^{-1})^*J_j$. (This 
function is supported near
$\gamma\bar x_j$.)

We will always take 
$\mu\in (0,\mu_0)$ where $\mu_0$ is sufficiently
small, so in 
particular
the supports of all functions $J_{j,\gamma}$ are 
disjoint.
Denote
$$
J_0=(1-\sum_{j,\gamma}J_{j,\gamma}^2)^{1/2}.
$$
Clearly, 
$J_0\in C^\infty(\widetilde M)$. Note that there exists
$c_0>0$ such 
that
$V\ge c_0\mu^{2\kappa}$ on $\supp J_0$.

functions $J_{j,\gamma}$
operators.

Denote by $\{\psi_{m,j}|\,m=1,2,\dots\}$ an orthonormal 
system of
eigenfunctions
of the operator $K_j$ in $L^2(\R^n)$ where 
the coordinates in a
neighborhood of the origin
$0\in \R^n$ are 
identified with the chosen coordinates near $\bar
x_j$, so that 
$0$
corresponds to $\bar x_j$. The corresponding eigenvalues will 
be
denoted $\lambda_{m,j}$.

Let us 
define
$\phi_{m,j}=J_j\psi_{m,j},$
extended by $0$ outside of a fixed 
small ball centered at $\bar x_j$. Then
$\phi_{m,j}\in 
C_0^\infty(\widetilde M)$ and it is supported near $\bar 
x_j$.

\begin{lemma}\label{L:almost-orthog} If $1/3<\kappa<1/2$, 
then
the functions $\phi_{m,j}$ satisfy the following 
``almost
orthogonality" 
relations
\begin{equation*}
(\phi_{m,j},\phi_{m',j'})=\delta_{j,j'}(\delta_{m,m'}+O(\mu^\kappa)),
\end{equation*}
\begin{equation*}
(H\phi_{m,j},\phi_{m',j'})=\delta_{j,j'}(\lambda_{m,j}\delta_{m,m'}+O(\mu^{3\kappa-1})),
\end{equation*}
where 
$j=1,\dots,N$, and $m$  belongs to a finite set.
\end{lemma}

The 
proof is the same as the proof of Lemma 2.3 in \cite{Sh}. 
The
unboundedness of
$A$ does not matter because only a finite number 
of points $\bar x_j$
is involved.

Not define 
$\phi_{m,j,\gamma}=T_\gamma \phi_{m,j}$. Then
$\phi_{m,j,\gamma}$ is 
supported near
$\gamma\bar 
x_j$.

\begin{lemma}\label{L:almost-orthog-gamma} If 
$1/3<\kappa<1/2$, then
the functions $\phi_{m,j}$ satisfy the 
following ``almost
orthogonality" 
relations
\begin{equation*}
(\phi_{m,j,\gamma},\phi_{m',j',\gamma'})=
\delta_{\gamma,\gamma'}\delta_{j,j'}(\delta_{m,m'}+O(\mu^\kappa)),
\end{equation*}
\begin{equation*}
(H\phi_{m,j,\gamma},\phi_{m',j',\gamma'})=
\delta_{\gamma,\gamma'}\delta_{j,j'}(\lambda_{m,j}\delta_{m,m'}+O(\mu^{3\kappa-1})),
\end{equation*}
where 
$j=1,\dots,N$, $\gamma\in\Gamma$, and $m$  belongs to a finite 
set.
\end{lemma}

\begin{proof}
The first relation is obvious because 
the operator $T_\gamma$ is
unitary and moves
supports by the action 
of $\gamma$. To prove the second estimate note
that it is obvious
if 
$\gamma\ne\gamma'$. If $\gamma=\gamma'$, then we 
get
$$
(H\phi_{m,j,\gamma},\phi_{m',j',\gamma})=(HT_\gamma\phi_{m,j},T_{\gamma}\phi_{m',j'})=
(T_\gamma^{-1}HT_\gamma\phi_{m,j},\phi_{m',j'})=(H\phi_{m,j},\phi_{m',j'}),
$$
and 
we can use the previous Lemma.
    \end{proof}

Now we will define 
two closed linear subspaces in
$\Phi^\F_\lambda\subset L^2(\F)$
and 
$\Phi_\lambda\subset L^2(\widetilde M)$ as 
follows:
$$
\Phi_\lambda^\F={\rm span}\{\phi_{m,j}|\,\lambda_{m,j}\le 
\lambda\},
$$
$$
\Phi_\lambda={\rm 
span}^c\{\phi_{m,j,\gamma}|\,\lambda_{m,j}\le \lambda\},
$$
where 
${\rm span}^c$ stands for closed linear span. 
Clearly
$\dim\Phi^\F_\lambda<\infty$,
and $\Phi_\lambda$ is a closed 
$(\Gamma,\sigma)$-invariant subspace
in $L^2(\widetilde 
M)$.

\begin{lemma}\label{L:dim}
$\hskip1.4in \dim 
\Phi_\lambda^\F=\dim_\Gamma 
\Phi_\lambda=N(\lambda;K).$
\end{lemma}

\begin{proof}
Clearly 
$\Phi^\F_\lambda=\oplus_{j=1}^N\Phi^\F_{\lambda,j}$, 
where
$\Phi^\F_{\lambda,j}$
is spanned by $\{\phi_{m,j}\}$ with fixed 
$j$.
But we have $\dim \Phi^\F_{\lambda,j}=N(\lambda;K_j)$.
Indeed, 
the eigenfunctions $\{\psi_{m,j}|m=1,\dots,N(\lambda;K_j)\}$
are 
linearly
independent and real analytic, so the corresponding 
$\phi_{m,j}$ are
also linearly independent
because 
$\phi_{m,j}=\psi_{m,j}$ near $\bar x_j$. It follows that
$\dim 
\Phi_\lambda^\F=N(\lambda;K)$.

For any fixed $j$ denote by 
$\{\tilde\phi_{m,j}|\,m=1,\dots,N(\lambda;K_j)\}$
the orthonormal 
system which is obtained from the 
system
$\{\phi_{m,j}|\,m=1,\dots,N(\lambda;K_j)\}$ by the 
Gram-Schmidt
orthogonalization process.
Then 
$\{\tilde\phi_{m,j}|\,\lambda_{m,j}\le\lambda\}$ is an
orthonormal 
basis in  $\Phi_\lambda^\F$,
and 
$\{T_\gamma\tilde\phi_{m,j}|\,\lambda_{m,j}\le\lambda,\gamma\in\Gamma\}$
is 
an orthonormal basis in $\Phi_\lambda$.  Note 
that
$T_\gamma\tilde\phi_{m,j}\in C_0^\infty(\widetilde M)$ and it 
is
supported near $\gamma\bar x_j$.

Denote by  $P_\lambda^\F$ and 
$P_\lambda$ the orthogonal projections on
$\Phi_\lambda^\F$ and 
$\Phi_\lambda$ respectively,
$K_\lambda^\F$ and $K_\lambda$ their 
Schwartz kernels. 
Then
$$
K_\lambda^\F(x,y)=
\sum_{\{m,j|\,\lambda_{m,j}\le\lambda\}}\tilde\phi_{m,j}\otimes\overline{\tilde\phi_{m,j}}=
\sum_{\{m,j|\,\lambda_{m,j}\le\lambda\}}\tilde\phi_{m,j}(x)\overline{\tilde\phi_{m,j}(y)}.
$$
In 
particular,
$$
K_\lambda^\F(x,x)=\sum_{\{m,j|\,\lambda_{m,j}\le\lambda\}}|\tilde\phi_{m,j}(x)|^2.
$$
Similarly 
we 
find
$$
K_\lambda(x,y)=
\sum_{\{m,j,\gamma|\,\lambda_{m,j}\le\lambda\}}T_\gamma\tilde\phi_{m,j}\otimes
\overline{T_\gamma\tilde\phi_{m,j}}
$$
and
$$
K_\lambda(x,x)=\sum_{\{m,j,\gamma|\,\lambda_{m,j}\le\lambda\}}
|\tilde\phi_{m,j}(\gamma 
x)|^2.
$$
It follows that $K_\lambda(x,x)=K_\lambda^\F(x,x)$ for all 
$x\in\F$, therefore
$$
\dim_\Gamma 
\Phi_\lambda=\Tr_\Gamma(P_\lambda)=\int_\F K_\lambda(x,x)dx=
\int_\F 
K_\lambda^\F(x,x)dx=\Tr(P_\lambda^\F)=\dim\Phi_\lambda^\F=N(\lambda;K),
$$
which 
proves the lemma.
\end{proof}

\begin{prop} For any $R>0$ and 
$\kappa\in (0,1/2)$ there exist
$\mu_0>0$ and $C>0$ such that
for any 
$\lambda\le R$ and any 
$\mu\in(0,\mu_0)$
\begin{equation}\label{E:below1}
N_\Gamma(\lambda+C\mu^\kappa;H)\ge 
N(\lambda;K).
\end{equation}
\end{prop}

\begin{proof}
Note first 
that if the estimate \eqref{E:below1} holds with some
$\kappa>0$, it 
holds also for
all smaller values of $\kappa$.
Now we should argue as 
in the proof of Lemma 2.6 in \cite{Sh} to conclude 
that
$\Phi_\lambda\in\Dom(H)$ and
$$
(Hf,f)\le 
(\lambda+C\mu^{3\kappa-1})(f,f), \quad f\in\Phi_\lambda.
$$
Applying 
the variational principle (Lemma \ref{L:var}) we conclude
that the 
estimate
\eqref{E:below1} holds with $3\kappa-1$ instead of $\kappa$. 
It
remains to notice that
$3\kappa -1$ takes all values from 
$(0,1/2)$ when $1/3<\kappa<1/2$.
\end{proof}

\begin{cor}
For any 
$R>0$ and $\kappa\in (0,1/2)$ there exist $\mu_0>0$ and $C>0$ such 
that
for any $\lambda\le R$ and any 
$\mu\in(0,\mu_0)$
\begin{equation}\label{E:below2}
N_\Gamma(\lambda;H)\ge 
N(\lambda-C\mu^\kappa;K).
\end{equation}
\end{cor}

\begin{rem}
If 
$H$ is flat near each of the points $\bar x_j$, $j=1,\dots,N$,
then a 
better
estimate is possible. Namely, for any $R>0$, $\kappa\in 
(0,1/2)$ and $\eps>0$
there exist $\mu_0>0$ and $C>0$
such that for 
any $\lambda\le R$ and any $\mu\in(0,\mu_0)$
\begin{equation}\label{E:below3}
N_\Gamma(\lambda;H)\ge 
N(\lambda-C\exp(-C^{-1}\mu^{-1+\eps});K).
\end{equation}
(Arguments 
in Sect.2 of \cite{Sh} apply here 
too.)

\end{rem}

\subsection{Estimate from above} Here we will prove 
an estimate from above
for $N_\Gamma(\lambda;H)$, similar to the 
estimate \eqref{E:below3}.

\begin{lemma}\label{L:above}
Assume that 
there exists an operator $D\in\U_{\widetilde 
M}(\Gamma,\bar\sigma)$
(i.e. a bounded  operator in $L^2(\widetilde 
M)$ commuting with all
magnetic translations
$T_\gamma$, 
$\gamma\in\Gamma$), such that

{\rm(a)} $\rank_\Gamma D\le \tilde 
k$;

{\rm(b)} $H+D\ge \tilde\lambda$.

\noindent
Then \ 
$N_\Gamma(\tilde\lambda-\eps;H)\le \tilde k$ for any 
$\eps>0$.
\end{lemma}

The proof does not differ from the proof of 
Lemma 3.7 in \cite{Sh}.

We will prove that for any $R>0$ there exist 
$C>0$ and $\mu_0>0$ such that
for any $\lambda<R$ and 
$\mu\in(0,\mu_0)$ there exists an operator
$D\in\U_{\widetilde 
M}(\Gamma,\bar\sigma)$, which satisfies the
conditions (a), (b)
with 
$\tilde k=N(\lambda;K)$, $\tilde\lambda=\lambda-C\mu^{1/5}$. It
would 
follow
from Lemma \ref{L:above} that
\begin{equation}\label{E:above1}
	N_{\Gamma}(\lambda- C\mu^{1/5};H)\le 	N(\lambda; K),
\end{equation}
hence
\begin{equation}\label{E:above}
	N_{\Gamma}(\lambda;H)\le 	N(\lambda+C\mu^{1/5}; K),
\end{equation}
which is the desired estimate from above.

Denote $E_\lambda^{(j)}(\mu)=\chi_{(-\infty,\lambda]}(K_j(\mu))$ which
is the spectral projection
of $K_j(\mu)$. It is an operator of the finite rank $N(\lambda;K_j)$
in $L^2(\R^n)$. Identifying
the coordinates in a neighborhood of 
$0\in\R^n$ with the chosen coordinates
near $\bar x_j$, so that 
$0\in\R^n$ corresponds to $\bar x_j$, we can
form an 
operator
$D_j=LJ_j E_\lambda^{(j)}(\mu) J_j$ in $L^2(\widetilde M)$, 
where
$L\in\R$, $L>\lambda$.
Clearly $D_j$ has a smooth Schwartz 
kernel supported in a
neighborhood of $(\bar x_j,\bar x_j)$
in 
$\widetilde M\times\widetilde M$.

Now denote
$$
\D_\F=\sum_{j=1}^n 
D_j
$$
and
$$
D=\sum_{\gamma\in\Gamma}T_\gamma D_\F 
T_\gamma^{-1}.
$$
It is easy to check that $T_\gamma D=D T_\gamma,\ 
\forall \gamma\in\Gamma$.

Note that $\rank 
E_\lambda^{(j)}(\mu)=N(\lambda;K_j)$, hence $\rank
D_j\le 
N(\lambda;K_j)$
and $\rank D_\F\le N(\lambda;K)$. (In fact, it is 
easy to see that we
have equalities in both
inequalities above, but 
we do not need this.) Also $D_\F$ has a
Schwartz kernel 
in
$C_0^\infty(\F\times \F)$.

Denote by $P_\F$ the orthogonal 
projection on the image  of $D_\F$. Choosing
an orthonormal basis in 
this image we see that the Schwartz kernel of $P_\F$
is also in 
$C_0^\infty(\F\times \F)$. For the orthogonal projection
$P$ on the 
closure
of the image of $D$ we 
have
$$
P=\sum_{\gamma\in\Gamma}T_\gamma P_\F 
T_\gamma^{-1}.
$$
Clearly $\rank D_\F=\Tr P_\F$ and $\rank_\Gamma 
D=\Tr_\Gamma P$.
But the  argument from the proof of Lemma 
\ref{L:dim} of the previous
subsection shows that $\Tr_\Gamma P=\Tr 
P_\F$, so
$$
\rank_\Gamma D=\rank D_\F\le N(\lambda;K),
$$
which 
proves the condition (a) with $\tilde k=N(\lambda;K)$.

To verify the 
condition (b) we will use the IMS localization
technique adopting 
the manifold version explained in \cite{Sh}.

Using the same 
functions $J_0,J_{j,\gamma}$ as in the previous section, we have
on 
$\widetilde M$:
$$
J_0^2+\sum_{j,\gamma}J_{j,\gamma}^2=1.
$$
We will also identify the functions $J_0,\ J_{j,\ga}$ with the
corresponding multiplication operators.

\begin{lemma}[The IMS localization formula] \label{identity}  The
following operator identity is true:
\begin{align*}
H &=J_0 H J_0+\sum_{j,\ga}J_{j,\ga}H J_{j,\ga}
+{\frac12}[J_0,[J_0,H]]+
{\frac12}\sum_{j,\ga}[J_{j,\ga},[J_{j,\ga},H]] \tag{L}\\
&=J_0 H J_0+\sum_{j,\ga}J_{j,\ga}H J_{j,\ga}
-\mu a^{(2)}(x,dJ_0(x))-\mu \sum_{j,\gamma}a^{(2)}(x,dJ_{j,\gamma}(x)),
\end{align*}
where $a^{(2)}$ is the principal symbol of $H$, considered as a
function on $T^*\widetilde M$.
\end{lemma}

The proof can be found e.g. in \cite{Sh}.

The last two terms in (L) are $\Gamma$-invariant and easily estimated
as $O(\mu^{1-2\kappa})$, so we have
\begin{equation}\label{E:est1}
H\ge J_0 H J_0+\sum_{j,\ga}J_{j,\ga}H J_{j,\ga}-C\mu^{1-2\kappa}I.
\end{equation}
Since $\mu H_A\ge 0$, we also have with some constant $c>0$
\begin{equation}\label{E:est2}
J_0 H J_0\ge \mu^{-1}VJ_0^2\ge c\mu^{-1+2\kappa} J_0^2,
\end{equation}
if $\mu\in(0,\mu_0)$, and $\mu_0$ is sufficiently small.
Note that the coefficient in the right hand side here
tends to $+\infty$ as $\mu\downarrow 0$.

Following Lemma 3.4 in \cite{Sh} we find that
\begin{equation}\label{E:est3}
J_jHJ_j\ge (1-C\mu^\kappa)J_jK_j(\mu)J_j-C\mu^{3\kappa-1}J_j^2\;.
\end{equation}
But we also have
\begin{align*}
J_{j,\gamma}HJ_{j,\gamma}=J_{j,\gamma}T_\gamma HT_\gamma^{-1}J_{j,\gamma}=
T_\gamma J_{j}HJ_{j}T_\gamma^{-1}\\
\ge  (1-C\mu^\kappa)T_\gamma
J_jK_j(\mu)J_jT_\gamma^{-1}-C\mu^{3\kappa-1}T_\gamma J_j^2
T_\gamma^{-1}\\
=(1-C\mu^\kappa)T_\gamma
J_jK_j(\mu)J_jT_\gamma^{-1}-C\mu^{3\kappa-1}J_{j,\gamma}^2\;.
\end{align*}
Let us sum over $j,\gamma$ and add $D$ and also $J_0 H J_0$. 
Using
the inequality
$$
K_j(\mu)+LE^{(j)}_\lambda\ge \lambda 
I,
$$
and also the estimates 
\eqref{E:est1},\eqref{E:est2},\eqref{E:est3}
above, we obtain 
then
\begin{equation*}
H+D\ge
c\mu^{-1+2\kappa}J_0^2+(1-C\mu^\kappa)\lambda\sum_{j,\gamma}J_{j,\gamma}^2-
C\mu^{3\kappa-1}\sum_{j,\gamma}J_{j,\gamma}^2-C\mu^{1-2\kappa}I.
\end{equation*}
Choosing 
here $\kappa=2/5$ we obtain
\begin{equation}\label{E:main-est}
H+D\ge 
(\lambda-C\mu^{1/5})I
\end{equation}
with a constant $C$ and with 
$\mu\in (0,\mu_0)$ for a sufficiently
small $\mu_0$.
This proves 
condition (b) and ends the proof of Theorem 2.1.
\hskip4.15in 
$\square$

\bigskip
\begin{rem} If $H$ is flat near all points $\bar 
x_j$, then the
estimate \eqref{E:above}
can be improved as follows: 
for any $R>0$ and $\eps>0$ there exist
$C>0$ and $\mu_0>0$
such that 
for all $\lambda<R$ and $\mu\in (0,\mu_0)$
$$
N_\Gamma(\lambda;H)\le 
N(\lambda+\mu^{1-\eps};K).
$$
Together with the improved estimate 
from below this provides the inclusion
$$
\spec(H)\cap 
(-\infty,R]\subset 
\bigcup_{p=1}^\infty
(\alpha_p-C\mu^{1-\eps},\alpha_p+C\exp(-C^{-1}\mu^{-1+\eps})).
$$
The 
necessary arguments can be found in 
\cite{Sh}.
\end{rem}

\end{document}